\documentclass[11pt]{article}
\usepackage{amsthm, amsmath, amssymb, amsfonts, url, booktabs, tikz, setspace, fancyhdr, bm}
\usepackage{hyperref}
\usepackage{geometry}
\usepackage{enumerate}
\usepackage[shortlabels]{enumitem}
\usepackage[babel]{microtype}
\usepackage[english]{babel}
\usepackage[capitalise]{cleveref}

\usepackage{tkz-graph,subcaption}
\usepackage{csquotes,mathrsfs}
\usetikzlibrary{patterns}
\usetikzlibrary{shapes}
\usetikzlibrary{decorations.pathreplacing}
\usetikzlibrary{decorations.pathmorphing}

\usepackage[T1]{fontenc}
\usepackage{amscd}
\usepackage{graphicx}
\usepackage{verbatim}
\usepackage{caption}
\usepackage{pdfpages,xcolor,framed,color}

\newcommand*{\myproofname}{Proof}

\newtheorem{thr}{Theorem}[section]
\newtheorem{ex}[thr]{Example}
\newtheorem{q}[thr]{Question}

\newtheorem{lem}[thr]{Lemma}
\newtheorem{prop}[thr]{Proposition}
\newtheorem{conj}[thr]{Conjecture}
\theoremstyle{definition}

\newtheorem*{defi*}{Definition}
\newtheorem{cor}[thr]{Corollary}

\newtheorem{remark}[thr]{Remark}

\newtheorem{claim}[thr]{Claim}

\def\A{\mathcal{A}}

\newcommand{\res}{\operatorname{res}}
\newcommand{\ext}{\operatorname{ext}}
\newcommand{\core}{\operatorname{core}}

\newcommand{\ny}{\mathcal N_{Y}}
\newcommand{\nx}{\mathcal N_{X}}
\newcommand{\diam}{diam}

\newcommand*{\abs}[1]{\lvert #1\rvert}

\def\a{\alpha}

\newcommand*{\av}{av}

\newenvironment{poc}{\begin{proof}[Proof of claim]}{\end{proof}}

\title{
On the extrema of the mean subtree order of graphs}

\author{
Stijn Cambie \thanks{Department of Computer Science, KU Leuven Campus Kulak-Kortrijk, 8500 Kortrijk, Belgium. Supported by FWO grants with grant numbers 1225224N and 1222524N. E-mail: \protect\href{mailto:stijn.cambie@hotmail.com}{\protect\nolinkurl{stijn.cambie@hotmail.com}} and \protect\href{mailto: 
                jorik.jooken@kuleuven.be}{\protect\nolinkurl{jorik.jooken@kuleuven.be}}}
\and
Jorik Jooken \footnotemark[1] \and Stephan Wagner\thanks{Institute of Discrete Mathematics, TU Graz, Austria, and Department of Mathematics, Uppsala University, Uppsala, Sweden. Supported by the Swedish research council (VR), grant 2022-04030. E-mail: \protect\href{mailto:stephan.wagner@tugraz.at}{\protect\nolinkurl{stephan.wagner@tugraz.at}}}
		}

\begin{document}

\maketitle

\begin{abstract}
 It has been conjectured that the minimum and maximum of the mean subtree order among connected graphs of order $n$ are attained by the path $P_n$ and clique $K_n$, respectively. Extending ideas due to Haslegrave and Vince, we confirm that the minimum is indeed attained by $P_n$. 
We also show that the maximum is attained by $K_n$ by proving that the ratio between spanning and almost spanning trees is maximised by the clique and applying a double-counting argument.
 % On the other hand, we discuss different approaches (both promising and flawed) that could lead to a proof of the extremality of $K_n$. 
\end{abstract}

\section{Introduction}

The mean subtree order, i.e., the average order of the subtrees of a tree, was first studied by Jamison about 40 years ago, see \cite{jamison_average_1983} and \cite{jamison_monotonicity_1984}.
The study of the mean subtree order (where a subtree of a graph is any subgraph that is a tree, i.e., not necessarily induced) was extended to arbitrary graphs by Chin, Gordon, MacPhee and Vincent~\cite{CGMV18}. The concluding section of~\cite{CGMV18} contains a number of open questions and conjectures, most of which have since been addressed, see~\cite{CM21, LXT24, Wagner21_prob_spanning}.
The main extremal question concerning the mean subtree order asks for the minimum and maximum values among all connected graphs of given order $n$. The answer was conjectured to be as stated below (see for example the introduction in~\cite{John22b}). In the following, we write $\mu(G)$ for the mean subtree order of a (connected) graph $G$.

\begin{conj}\label{conj:main}
 Let $G$ be a connected graph of order $n$. Then 
 \begin{enumerate}
 \item\label{itm:min_Pn} $\mu(P_n) \le \mu(G)$ with equality if and only if $G=P_n$,
 \item\label{itm:max_Kn} $\mu(G) \le \mu(K_n)$ with equality if and only if $G=K_n$.
 \end{enumerate}
\end{conj}

Let us remark that the restricted version of~\ref{itm:min_Pn}, where $G$ is required to be a tree, was already proven by Jamison~\cite{jamison_average_1983}. We combine ideas from~\cite{John22, Vince22} to prove~\ref{itm:min_Pn} in~\cref{sec:proofPnMinimiser}.

\begin{thr}\label{thr:Pn_istheminimiser}
    Let $G$ be a connected graph of order $n$.
    Then $\mu(G) \ge \frac{n+2}{3}$ with equality if and only if $G=P_n.$
\end{thr}

Moreover, we prove~\cref{conj:main}~\ref{itm:max_Kn} through the following auxiliary result on the ratio of spanning trees and subtrees of order $n-1$. Let $s_k(G)$ denote the number of (not necessarily induced) subtrees of $G$ with $k$ vertices.

\begin{prop}\label{redconj_Kn}
 For a graph $G$ of order $n\ge 2$, we have 
$$s_{n-1}(G) \ge \frac{s_{n-1}(K_n)}{s_n(K_n)} s_n(G).$$
\end{prop}

In~\cref{sec:maximum_Kn}, we prove~\cref{redconj_Kn}, and show that it implies~\cref{conj:main}~\ref{itm:max_Kn}, and also that $K_n$ maximises $p(G)$, the fraction of subtrees that are spanning (\cref{conj2}).

\begin{thr}\label{thr:Kn_isthemaximiser}
    Let $G$ be a connected graph of order $n$.
    Then $\mu(G) \le \mu(K_n)$ with equality if and only if $G=K_n.$
\end{thr}

We further consider some conjectures that would imply~\cref{conj:main}~\ref{itm:min_Pn} or~\ref{itm:max_Kn}, or both, and give an overview of some directions that are of independent interest.

A potential approach would be based on adding or removing edges. Let $G+e$ denote the graph $G$ with the edge $e$ added. More generally, for a set of edges $M$, the graph $G+M$ is the graph with vertex set $V(G)$ and edge set $E(G) \cup M$. Chin, Gordon, MacPhee and Vincent conjectured that for every two non-adjacent vertices $u,v$ in a graph $G$, adding the edge $e=uv$ increases the mean subtree order, i.e., $\mu(G+e)>\mu(G)$ (see~\cite[Conj.~7.4]{CGMV18}). The upper bound (\cref{conj:main}~\ref{itm:max_Kn}) would be an easy corollary, and the lower bound (\cref{conj:main}~\ref{itm:min_Pn}) would follow from Jamison's result that the path is extremal among trees. However, the conjecture was disproved by Cameron and Mol~\cite{CM21}, who provided explicit counterexamples.

Recently, two groups independently generalised these counterexamples by proving that for every positive integer $k$, there is a graph and a choice of $k$ edges whose addition decreases the mean subtree order~\cite{CCHT23,CWL23}.\footnote{Here~\cite{CCHT23} was submitted first and~\cite{CWL23} was published first.}
A conjecture presented in the conclusion of~\cite{CWL23} is refuted by one of the examples in~\cite{CCHT23}, where it is demonstrated that the mean subtree order can decrease by as much as $\frac{1}{3}(1-o(1))n$ when $k$ edges are added. We provide an even more robust counterexample to the conjecture presented in~\cite[Conj.~7.4]{CGMV18} by proving the following theorem in~\cref{sec:nonuniversal}. 

\begin{thr}\label{thr:strong_counterexamples_toConj7.4Chin}
 There are graphs $G$ of order $n$ and size $o(n^2)$ for which there is only one edge $e$ up to isomorphism ($o(n)$ many out of the $\Theta(n^2)$ edges in $G^c$) for which $\mu(G+e)>\mu(G)$. 
 There are also graphs for which adding any maximal matching $M$ in $G^c$ reduces the mean subtree order, i.e., $\mu(G+M)<\mu(G)$.
\end{thr}

\cref{conj:main} would also be a corollary of the following conjecture\footnote{Note that~\cref{conj:main2}\ref{itm:main2_Kn} has been stated before in~\cite[Conj.~1.2]{CM21}.}. Here, $G \setminus e$ denotes the graph obtained by deleting the edge $e$ in the graph $G$.

\begin{conj}\label{conj:main2}
 Let $G$ be a connected graph of order $n$. 
 \begin{enumerate}
 \item\label{itm:main2_Pn} If $G$ is not a tree, then there exists an edge $e \in G$ such that $G \setminus e$ is connected and $\mu(G \setminus e)<\mu(G).$
 \item\label{itm:main2_Kn} If $G \neq K_n$, then there exist two non-adjacent vertices $u,v$ such that $\mu(G+uv)>\mu(G).$
 \end{enumerate}
\end{conj}

In~\cref{sec:nonuniversal} we show that both statements are false when the existence is replaced by a ``for all'' statement.

Another direction that would imply~\cref{conj:main}~\ref{itm:min_Pn} would be a two-step argument in which one proves some monotonicity results analogous to those obtained in~\cite{jamison_monotonicity_1984} for trees. 
However, this approach also fails.
In~\cref{sec:nomonotonicity}, we show that analogous monotonicity results do not hold for arbitrary connected graphs.

Another interesting direction is related to Jamison's contraction conjecture. We let $G/e$ denote the graph $G$ with the edge $e$ contracted.
\begin{conj}\cite[Conj.~5.6]{jamison_monotonicity_1984}\label{conj:contractionTrees}
 In a tree $T$, the contraction of any edge $e \in E(T)$ reduces the mean subtree order by at least $\frac 13$; $\mu(T)-\mu(T /e) \ge \frac 13$.
\end{conj}

This conjecture was proven in~\cite{LXWW23, Ruoyu24} (with equality if and only if $T$ is a path). It is tempting to generalise the conjecture to arbitrary connected graphs.

\begin{conj}\label{conj:contraction}
 In a connected graph $G$, the contraction of any edge $e \in E(G)$ reduces the mean subtree order by at least $\frac 13$; $\mu(G)-\mu(G /e) \ge \frac 13,$ with equality if and only if $G$ is a path.
\end{conj}

Proving~\cref{conj:contraction} would immediately lead to an alternative proof of~\cref{conj:main}~\ref{itm:min_Pn}. As we will see in~\cref{sec:nomonotonicity}, it is also sufficient to prove~\cref{conj:contraction} in the case that $e$ is a bridge of the graph $G.$

We close this section by briefly discussing another generalisation of the study of the mean subtree order of trees. For a graph $G=(V,E)$, we call a vertex set $S \subseteq V$ connected if $G[S]$ is connected and denote by $\av(G)$ the average size of a connected vertex set of a graph. Note that for a tree $T$, we have $\mu(T)=\av(T)$. As mentioned before, Jamison~\cite{jamison_average_1983} proved that the path $P_n$ attains the minimum mean subtree order among all trees, and it is also known that the path minimises $\av(G)$ among all connected graphs, as proven independently by Haslegrave and Vince~\cite{John22,Vince22}. However, the situation differs for the maximum of $\av(G)$ because the clique $K_n$ is not extremal for sufficiently large $n$ (since for example, $\frac{n}{2} \sim \av(K_n) < \av(S_n) \sim \frac{n+1}{2}$ for sufficiently large $n$). This raises the question of whether the extremal graph (attaining the largest value of $\av(G)$ among graphs $G$ of order $n$) must be a tree, which, according to the caterpillar conjecture~\cite[Prob.~(7.1)]{jamison_average_1983}, would be a caterpillar.

We also mention a fundamental question (alluded to in~\cite{John22b,John22}) that connects $\mu(G)$ and $\av(G)$.

\begin{q}\label{mu_vs_av}
    Is it true that $\mu(G) \ge \av(G)$ holds for every connected graph $G$?
\end{q}

If true, this would imply~\cref{conj:main}~\ref{itm:min_Pn} immediately from the aforementioned result of Haslegrave and Vince~\cite{John22, Vince22}.

%We have proven one part of the main conjecture, namely~\cref{conj:main}~\ref{itm:min_Pn}, leaving~\cref{conj:main}~\ref{itm:max_Kn} open.
%The alternative strategies in~\cref{conj:main2},~\cref{conj:contraction} (for bridges),~\cref{redconj_Kn} (and~\cref{conj2}) and~\cref{mu_vs_av} may be interesting in their own right. We also draw the reader's attention to a related open problem in~\cite[Prob.~5.1]{CM21}, which asks about the minimum and maximum change in the density of a graph due to edge addition.
\section{The path minimises the mean subtree order among connected graphs}\label{sec:proofPnMinimiser}

Let $N(G)=\sum_{k=1}^n s_k(G)$ be the number of subtrees of $G$ and $R(G)= \sum_{k=1}^n k\cdot s_k(G) $ be the total order of all subtrees of $G$.
Similarly, $N(G,T)$ and $R(G,T)$ denote the number and total order of subtrees of $G$ containing $T.$

For a graph $G$ of order $n$, we let $\mu(G) = \frac{R(G)}{N(G)} = \frac{ \sum_{k=1}^n k\cdot s_k(G)}{ \sum_{k=1}^n s_k(G)}$ be the mean subtree order of the graph $G$. 
The local mean subtree order in a vertex $v$ or subtree $T$, $\mu(G, v)$ or $\mu(G, T)$, is the average order of all subtrees containing $v$ resp.~$T$.

We prove the following variant of~\cite[Thm.~3.1]{Vince22}, whose proof is along the same lines.
 
\begin{prop}\label{prop:localsubtreeorder} If $T$ is a (non-empty) subtree of a connected graph $G$ of order $n$, then
\[\mu(G,T) \ge \frac{n+\abs{T}}{2}.\]
\end{prop}

\begin{proof} 
The proof is by induction on the integer $d = n- \abs{T}$. If $d=0$, then $\mu(G,T)=n=\frac{n+n}{2}$, so the base case of the induction is true. 
Assume that the statement is true for $d-1$, and let $(G,T)$ be a pair of a graph and a subtree for which $n-\abs T = d$. 
Let $Q= \bigcup_{v \in T} N(v) \setminus V(T)$ be the set of vertices that are adjacent to a vertex in $T$ but are not in $T$. 
We take a vertex $x$ in $Q$ and define $\mathcal T$ as the set of all minimal subtrees (with respect to inclusion) that contain both $T$ and $x.$ We let $y$ be a (fixed) neighbour of $x$ that belongs to $T$.

Let $G' = G \setminus x$, let $G_1$ be the connected component of $G \setminus x$ containing $T$, and let $G_2$ be the union of the other components. Then
\[N(G,T)= \sum_{S \in \mathcal T}  N(G, S )+ N(G',T) =\sum_{S \in \mathcal T} N(G, S ) + N(G_1,T)\]
and
\[R(G,T)= \sum_{S \in \mathcal T}  R(G, S )+ R(G',T) =\sum_{S \in \mathcal T} R(G, S ) + R(G_1,T).\]
Let us write $m$ for the number of vertices of $G_2$ %Denote the vertex set of $G_2$ by $V_2$, and let $m = \abs{V_2}$
(which is $0$ if $x$ is not a cut vertex). Observe that
$$\sum_{S \in \mathcal T} N(G, S ) \ge N(G, T \cup \{xy\}) \ge (m+1) N(G_1\cup \{xy\}, T \cup \{xy\}) = (m+1) N(G_1, T).$$
Here, the second inequality is true since there are at least $m+1$ ways to extend a subtree of $G_1\cup \{xy\}$ to a subtree of $G$.
For this, consider a spanning tree of $G \setminus G_1$ and a sequence of subtrees starting with $x$ and ending with the spanning tree, obtained by adding one edge at a time.

Now, using the induction hypothesis, we get
\begin{align*} R(G,T) &= \sum_{S \in \mathcal T} R(G, S ) + R(G_1,T)\\
&\ge \frac{n+(\abs{T}+1)}{2} \sum_{S \in \mathcal T} N(G, S ) + \frac{(n-1-m)+\abs T}{2} N(G_1,T) \\&=
\frac{n+\abs{T}}{2} \left( \sum_{S \in \mathcal T} N(G, S ) + N(G_1,T)  \right ) + \frac{1}{2} \Big (\sum_{S \in \mathcal T} N(G, S )- (m+1) N(G_1, T) \Big ) \\ & \ge \frac{n+\abs{T}}{2} N(G, T),
\end{align*}
completing the induction and thus the proof.
\end{proof} 

Next, we prove the analogue of~\cite[Lem.~3]{John22}, again by modifying the original proof to the setting of the mean subtree order.
Note that $N(G)$ denotes the number of subtrees in our paper, while in~\cite{John22}, it denotes the number of connected induced subsets.
The main adaptation lies in the analogue of~\cite[Clm.~3.3]{John22}, which is our~\cref{clm:ell_ny}.

\begin{prop}\label{Fraction}
 Let $G$ be any connected graph on $n \ge 3$ vertices. Then G contains a vertex v such that $G \setminus v$ is connected and $N(G, v) \ge 2\frac{N(G)}{n+1}$, with equality if and only if $G$ is a path.
\end{prop}

\begin{proof}
    Let $P=v_0v_1 \cdots v_{d}$ be a diameter of $G$ (i.e., $d=\diam(G)$, where we can assume that $2 \le d \le n-2$ since for the path and complete graph, the statement is easily seen to be true).
    Let $V_i(P)$ denote the set of vertices of $V(G) \setminus V(P)$ that are at distance $i$ from $P$, and let its elements be denoted by $v_i^j$ for $1 \le j \le \abs{V_i(P)}$. The specific order of the vertices in $V_i(P)$ is irrelevant as long as it is kept consistent throughout the proof.
    We now define an auxiliary coloured directed multidigraph $H$ whose vertex set is the set of subtrees of $G$.
    Let $S$ be a subtree of $G$. If $d(S,P)=i\ge 2$ then choose a vertex $v_{i-1}^j$ which is at distance $1$ from $S$, with minimal index $j$.
    Let $v_i^{k}$ be a neighbour of $v_{i-1}^j$ belonging to $S$, with minimal index $k$.
    Let $S'=S \cup \{v_{i-1}^jv_i^{k}\}$ be the tree $S$ with the edge $v_{i-1}^jv_i^{k}$ added.
  Add two directed edges, one red and one blue, from $S$ to $S'$.
  If $d(S,P)\le 1$ and $v_0\not\in S$, let $i$ be minimal such that $d(v_i,S)=1$. Add a red edge from $S$ to $S\cup\{v_iv_{i+1}\}$ if $v_{i+1} \in S$, and otherwise to the subtree $S\cup\{v_iv_1^j\}$ with minimal $j$. 
  If $d(S,P)\le 1$ and $v_d\not\in S$, let $j$ be maximal such that $d(v_j,S)=1$. Add a blue edge from $S$ to $S\cup\{v_{j-1}v_j\}$ if $v_{j-1} \in S$, and otherwise to the subtree $S\cup\{v_1^kv_j\}$ with maximal $k$.
  
  For every subtree $S$ and colour $c \in \{\text{red, blue}\},$ there is at most one outgoing and at most one incoming edge of colour $c$.
  The subdigraph of $H$ induced by all edges of colour $c$ is thus a disjoint union of directed paths.
    A vertex of $H$ (a subtree of $G$) with no red (resp.~blue) incoming edge is called a red (resp.~blue) top.
    We write $\mathcal T$ for the set of all pairs $(S,c)$ where $S$ is a subtree of $G$ with no incoming edge in colour $c$. For $\tau=(S,c)\in\mathcal T$ we write $\ell(\tau)$ for the length of the path of colour $c$ from the vertex corresponding to $S$ in $H$.
    For a subcollection $\mathcal T'\subset \mathcal T$ of tops with associated colours, we use $\ell(\mathcal T')=\frac{1}{\abs{ \mathcal T'} } \sum_{\tau \in \mathcal T'} \ell(\tau)$ to denote the average length of the paths corresponding to the tops in $\mathcal T'$.

    We divide $\mathcal T$ into three parts. We say that a top $(S,c)$ is ``high'' if $V(S) \subseteq V(G)\setminus V(P)$, ``low'' if $V(S) \subseteq V(P)$, and ``normal'' otherwise (i.e., if $S$ contains vertices from both $V(P)$ and $V(G)\setminus V(P)$). We write $\mathcal{H}$, $\mathcal{L}$ and $\mathcal{N}$ for the sets of high, low and normal tops respectively. 
    For a subtree $S$ of $G$, we define its residue $\res(S)$ as the forest $S[V(G) \setminus V(P)]=S \setminus V(P).$
    
    For a forest $X$ with $d(X,V(P))=1$, let $i_X$ be the minimal index $i$ such that $v_i$ has a neighbour in $X$, and let $j_X$ be the maximal such index. 
    Next, we define the core---denoted $\core(S)$---of a subtree $S$ (which will be a subtree of $S$) in two steps.
    First, take $S'=S \setminus \{v_k \mid k < i_{\res(S)} \vee k > j_{\res(S)} \}$, i.e., the subtree of $S$ spanned by the vertices different from vertices $v_k$ for which $k < i_{\res(S)}$ or $ k > j_{\res(S)}.$
    Second, if $v_{i_{\res(S)}}$ or $v_{j_{\res(S)}}$ is a leaf of $S'$, we remove it from $S'$, otherwise it is retained. The final result obtained from $S'$ after possibly removing one or both of these vertices is $\core(S)$.

    For a normal top $\tau=(S,c),$ we also write $\res(\tau)$ and $\core(\tau)$ for $\res(S)$ and $\core(S)$ respectively.
    For a subtree $Y$ of $G$, define $\ny=\{ \tau \in \mathcal N \colon \core(\tau)=Y\}$ as the collection of all normal tops whose core equals the tree $Y$.
    Let $\ny^r$ and $\ny^b$ be the partition into the red and blue tops of $\ny$, respectively.
    Note that $\res(\tau)=\res(\core(\tau)).$ 

    \begin{claim}\label{clm:ell_ny}
        Assume that $\ny \neq \emptyset$. Then $\ell(\ny) \le \frac{d}{2} $ if $i_{\res(Y)}=j_{\res(Y)}$, and $\ell(\ny) \le \frac{d+1}{2} $ if $i_{\res(Y)}<j_{\res(Y)}.$
    \end{claim}
    \begin{poc}
        Let $X=\res(Y).$
        First assume that $i_X=j_X.$
        We distinguish two cases, based on whether the core contains $v_{i_X}$ or not.
        \begin{itemize}
            \item If $\deg_Y(v_{i_X})\ge 2$ and thus $v_{i_X} \in Y$, the red tops in $\ny$ are exactly equal to $Y$ extended by a path $v_{i_X}\cdots v_b$ for some $i_X \le b \le d,$ and the blue tops in $\ny$ are exactly equal to $Y$ extended by a path $v_a \cdots v_{i_X}$ for some $0 \le a \le i_X.$
            Thus there are $i_X+1$ blue tops $\tau$ with length $\ell(\tau) = d-i_X$ and $d-i_X+1$ red tops $\tau$ with length $\ell(\tau) = i_X.$
            As such, by the elementary inequality $(d-2i_X)^2 \ge 0$, which is equivalent to $4i_X(d-i_X)\le d^2$, we have
            $$\ell(\ny)=\frac{ (i_X+1)(d-i_X)+(d-i_X+1)i_X}{i_X+1+d-i_X+1} = \frac{2i_X(d-i_X) +d}{d+2} \le \frac{ d}{2}. $$

            \item Otherwise, $Y=X$. Let $x$ be the number of neighbours of $v_{i_X}$ in $X$.
            The number of red tops in $\ny$ is equal to $x(d-i_X+1)-1,$ i.e., $\abs{\ny^r}=x(d-i_X+1)-1.$
            This is because adding any edge between $v_{i_X}$ and a neighbour in $X$, together with some path starting in $v_{i_X}$ towards $v_{d}$, results in a red top, except for the red outneighbour of $X$ in $H$ ($X$ with the edge from $v_{i_X}$ to the neighbour $v_1^k$ with minimal index $k$ added).
            Analogously, $\abs{\ny^b}=x(i_X+1)-1.$ %the number of blue tops in $\ny$
            So in this case, using the same inequality as before,
            $$\ell(\ny)=\frac{ [x(i_X+1)-1](d-i_X)+[x(d-i_X+1)-1]i_X}{x(d+2)-2} = \frac{2i_X(d-i_X)x +d(x-1)}{xd+2(x-1)} \le \frac{ d}{2}. $$   
        \end{itemize}
        Next, we assume $i_X<j_X.$
        Now we have to consider three cases, which are depicted in~\cref{fig:3sit}.
        \begin{itemize}
            \item $v_{i_X}, v_{j_X} \not \in Y$\\
            Let $x$ and $y$ be the number of neighbours of $v_{i_X}$ resp.~$v_{j_X}$ in $Y$, and let $z$ be the number of common neighbours of these two vertices.

            First, we assume that $Y=X$, i.e., there is no vertex of $P$ in $Y$.
            
            The number of blue normal tops with core $Y$ can be counted by considering all the possible extensions of $Y$. There are 
            \begin{itemize}
                \item $y-1$ choices to add an edge containing $v_{j_X}$ (length of top is $d-j_X$),
                \item $x(i_X+1)$ choices to add an edge between $Y$ and $v_{i_X}$ and possibly extend by a path towards $v_0$ (length of top is $d-j_X+1$),
                \item the above $x(i_X+1)$ possibilities, but for each of them additionally one of the $y-1$ choices to add an edge containing $v_{j_X}$ (length of top is $d-j_X$).
            \end{itemize}
            We conclude that $\abs{\ny^b}=xy(i_X+1)+y-1$ and analogously $\abs{\ny^r}=xy(d-j_X+1)+x-1$.
            Now
            $$\ell(\ny) = \frac{(xy(i_X+1)+y-1)(d-j_X)+x(i_X+1)+ (xy(d-j_X+1)+x-1)i_X+y(d-j_X+1) }{xy(i_X+1)+y-1+xy(d-j_X+1)+x-1}.$$
%             \begin{align*}
%                 \ell(\ny)&=\frac{\abs{\ny^b}(d-j_X)+x(i_X+1)+ \abs{\ny^r}i_X+y(d-j_X+1)}{xy(i_X+1)+y-1+xy(d-j_X+1)+x-1} \\
%                 &=\frac{(xy(i_X+1)+y-1)(d-j_X)+x(i_X+1)+ (xy(d-j_X+1)+x-1)i_X+y(d-j_X+1) }{xy(i_X+1)+y-1+xy(d-j_X+1)+x-1}.% \\
% %                &=\frac{2xy(i_X+1)(d-j_X+1)-x(y-2)(i_X+1)-(x-2)y(d-j_X+1)-i_X-(d-j_X)-x-y}{xy(i_X+1)+xy(d-j_X+1)+x+y-2}. 
%             \end{align*}
            Setting $a = i_X+1$ and $b = d-j_X+1$, this becomes
        \begin{align*}
        \ell(\ny) &= \frac{xy(2ab-a-b)+(2a-1)x+(2b-1)y-a-b+2}{xy(a+b)+x+y-2} \\
        &= \frac{a+b}{2} - \frac{xy(a-b)^2 + (2xy-x-y)(a+b)-2(x-y)(a-b)+(2x+2y-4)}{2(xy(a+b)+x+y-2)}.
        \end{align*}
If we can show that the numerator of the second fraction is nonnegative, the claim follows since $a+b = d+1 - (j_X-i_X-1) \le d+1$. To this end, observe that $(a-b)^2 \ge \abs{a-b}$ since $a$ and $b$ are integers, that $a+b \ge \abs{a-b}$, and that $2x+2y-4 \ge 0$. Thus the numerator is bounded below by
$$xy\abs{a-b} + (2xy-x-y)\abs{a-b} - 2\abs{x-y}\abs{a-b} = (3xy-x-y-2\abs{x-y})\abs{a-b}.$$
This is nonnegative since
$$3xy-x-y-2\abs{x-y} = 3xy - 3\max\{x,y\} +\min\{x,y\} \ge
3(xy - \max\{x,y\}) \ge 0.$$
%             If $x=y=1$, we have
%             $$\ell(\ny) = \frac{2(i_X+1)(d-j_X+1)}{(i_X+1)+(d-j_X+1)} \le \frac{(i_X+1)+(d-j_X+1)}2
%                 \le \frac{d+1}2$$
%             by the inequality between harmonic and arithmetic mean. If $x,y>1$, then the numerator of $\ell(\ny)$ is clearly less than $2xy(i_X+1)(d-j_X+1)$ while the denominator is greater than $xy(i_X+1)+xy(d-j_X+1)$, thus
% $$\ell(\ny) \le \frac{2(i_X+1)(d-j_X+1)}{(i_X+1)+(d-j_X+1)}\le \frac{(i_X+1)+(d-j_X+1)}2
%                 \le \frac{d+1}2$$
% once again. Finally, if $x=1$ and $y>1$ (the case $y=1, x>1$ is analogous), we use the substitution $a=i_X+1$ and $b=d-j_X+1$

% , we have that $\ell(\ny)\le \frac{2ab+b}{2(a+b)+1} \le \frac{a+b}{2}=\frac{d+1}2$. The central inequality is equivalent to $(b-a)^2 \ge \frac{b-a}{2}$, which is true for integers.

% we have that it is again bounded by $\frac{(i_X+1)+(d-j_X+1)}2\le \frac{d+1}2.$
% After a substitution, this is true by $\frac{2ab+b}{2(a+b)+1} \le \frac{a+b}{2}$ which is equivalent with $(b-a)^2 \ge \frac{b-a}{2} $, which is true for integers.

If $Y\neq X,$ both $(Y,\text{blue})$ and $(Y,\text{red})$ are tops as well.
But these two tops have length $(i_X+1)$ resp. $(d-j_X+1)$ with an average of $\frac{d-j_X + i_X + 2}{2} \leq \frac{d+1}{2}$, so the same upper bound is still valid.

            \item $v_{i_X} \in Y, v_{j_X} \not\in Y$ (the case $v_{j_X} \in Y, v_{i_X} \not \in Y$ is analogous)\\
            Let $y$ be the number of neighbours of $v_{j_X}$ in $Y$.
            Now, by similar reasoning as in the previous cases, $\abs{\ny^b}=y(i_X+1)$ and $\abs{\ny^r}=y(d-j_X+1)+1$.
            There are $(y-1)(i_X+1)$ blue tops that contain $v_{j_X}$ (length $d-j_X$) and $i_X+1$ that do not (length $d-j_X+1$).
            On the other hand, $Y$ itself is a red top, and all the others contain $v_{j_X}$. All the red tops have length $i_X$. 
            This implies that %for $y \ge 2$
   $$\ell(\ny)=\frac{y(i_X+1)(d-j_X)+(i_X+1)+y(d-j_X+1)i_X+i_X }{y(i_X+1)+y(d-j_X+1)+1}.$$
   Setting $a = i_X+1$ and $b = d-j_X+1$ again, this becomes
   \begin{align*}
        \ell(\ny) &= \frac{y a (b - 1) + a + y b (a - 1) + a-1}{ya+yb+1} \\
        &= \frac{a+b}{2} - \frac{y(a-b)^2 + 2y(a+b) + b - 3a+2}{2(ya+yb+1)}.   
   \end{align*}
   We can conclude as in the previous case, since the numerator of the second fraction is again nonnegative:
$$y(a-b)^2 + 2y(a+b) + b - 3a + 2 > |a-b| + 2(a+b) - b - 3a = |a-b| + (b-a) \geq 0.$$
   
            % \begin{align*}
            %     \ell(\ny)&=\frac{y(i_X+1)(d-j_X)+(i_X+1)+y(d-j_X+1)i_X+i_X }{y(i_X+1)+y(d-j_X+1)+1}\\
            %     &\le \frac{y(i_X+1)(d-j_X+1)+y(d-j_X+1)(i_X+1) }{y(i_X+1)+y(d-j_X+1)}\\
            %     &= \frac{2(i_X+1)(d-j_X+1) }{(i_X+1)+(d-j_X+1)}\\
            %     &\le \frac{(i_X+1)+(d-j_X+1)}2\\
            %     &\le \frac{d+1}2.
            % \end{align*}
            % Here we have used the inequality between the arithmetic mean and the harmonic mean (applied to $i_X+1$ and $d-j_X+1$).
            % For $y=1$, we use the substitution $a=i_X+1$ and $b=d-j_X+1$ and obtain
            % $$\ell(\ny) =\frac{2ab+a-b-1}{a+b+1}\le \frac{a+b}{2} \text{ since } (a-b-1)^2+a+b+1 \ge 0$$       
        \item $v_{i_X}, v_{j_X} \in Y$ (thus $\deg_Y(v_{i_X}), \deg_Y(v_{j_X})\ge 2$)\\
        This is the easiest case, where we can immediately find that $\abs{\ny^b}=i_X+1 = a$, $\abs{\ny^r}=d-j_X+1 = b$ and 
        $$\ell(\ny)=\frac{(i_X+1)(d-j_X)+(d-j_X+1)i_X}{i_X+1+d-j_X+1} 
        =\frac{a(b-1)+b(a-1)}{a+b} = \frac{a+b}{2} - \frac{(a-b)^2}{2(a+b)} - 1,$$
        and we are done as in the other two cases.\qedhere
        \end{itemize}
    \end{poc}

\begin{figure}[ht]
    \centering
    \begin{tikzpicture}[scale=0.6]

        \draw[fill=black!10!white, draw=black!10!white] (-0.5,4) -- (-0.5,1) -- (2,1) -- (2,-0.5) --
(4,-0.5) -- (4,1) -- (6.5,1)  -- (6.5,4) -- cycle;

        \foreach \x in {0,1,2,3,4}{
         \draw[fill] (\x,2) circle (0.1);
        \draw[dotted] (1,0)--(\x,2);
        }
        \foreach \x in {3,4,5,6}{
         \draw[fill] (\x,2) circle (0.1);
        \draw[dotted] (5,0)--(\x,2);
        }
\draw [decorate,decoration={brace,amplitude=4pt},xshift=0pt,yshift=0pt] (-0.2,2.12) -- (2.2,2.12) node [black,midway,yshift=0.35cm]{$x-z$};
\draw [decorate,decoration={brace,amplitude=4pt},xshift=0pt,yshift=0pt] (2.8,2.12) -- (4.2,2.12) node [black,midway,yshift=0.35cm]{$z$};
\draw [decorate,decoration={brace,amplitude=4pt},xshift=0pt,yshift=0pt] (4.8,2.12) -- (6.2,2.12) node [black,midway,yshift=0.35cm]{$y-z$};

     \draw[dashed] (-1,0)--(7,0);

        \draw[fill] (1,0) circle (0.1);
        \draw[fill] (5,0) circle (0.1);
        \draw[fill] (-1,0) circle (0.1);
        \draw[fill] (7,0) circle (0.1);

    \node at (-1,-0.38) {$v_0$};
        \node at (1,-0.38) {$v_{i_X}$};
         \node at (5,-0.38) {$v_{j_X}$};
          \node at (7,-0.38) {$v_{d}$};
          \node at (3,4.38) {$Y$};
    \end{tikzpicture}\quad
    \begin{tikzpicture}[scale=0.6]

        \draw[fill=black!10!white, draw=black!10!white] (-0.5,3.5) -- (-0.5,1) -- (1,1) -- (1,-1.1) --
(4,-1.1) -- (4,1) -- (5.5,1)  -- (5.5,3.5) -- cycle;

        \foreach \x in {3,4,5}{
         \draw[fill] (\x,2) circle (0.1);
        \draw[dotted] (5,0)--(\x,2);
        }

\draw [decorate,decoration={brace,amplitude=4pt},xshift=0pt,yshift=0pt] (2.8,2.12) -- (5.2,2.12) node [black,midway,yshift=0.35cm]{$y$};

        \foreach \x in {1,2}{
        
        \draw[dashed] (1.5,0)--(\x,1.5);
        }

     \draw[dashed] (-1,0)--(7,0);

        \draw[fill] (1.5,0) circle (0.1);
        \draw[fill] (5,0) circle (0.1);
        \draw[fill] (-1,0) circle (0.1);
        \draw[fill] (7,0) circle (0.1);

    \node at (-1,-0.38) {$v_0$};
        \node at (1.5,-0.38) {$v_{i_X}$};
         \node at (5,-0.38) {$v_{j_X}$};
          \node at (7,-0.38) {$v_{d}$};
          \node at (3,3.88) {$Y$};
    \end{tikzpicture}
    \quad
    \begin{tikzpicture}[scale=0.6]

        \draw[fill=black!10!white, draw=black!10!white] (-0.5,3.5) -- (-0.5,1) -- (0.5,1) -- (0.5,-1) --
(3.5,-1) -- (3.5,1) -- (4.5,1)  -- (4.5,3.5) -- cycle;

        \foreach \x in {0.75,1.5}{
        
        \draw[dashed] (1,0)--(\x,1.5);
        }
        \foreach \x in {2.5,3.25}{
        \draw[dashed] (3,0)--(\x,1.5);
        }

     \draw[dashed] (-1,0)--(5,0);

        \draw[fill] (1,0) circle (0.1);
        \draw[fill] (5,0) circle (0.1);
        \draw[fill] (-1,0) circle (0.1);
        \draw[fill] (3,0) circle (0.1);

    \node at (-1,-0.38) {$v_0$};
        \node at (1,-0.38) {$v_{i_X}$};
         \node at (3,-0.38) {$v_{j_X}$};
          \node at (5,-0.38) {$v_{d}$};
          \node at (2,3.88) {$Y$};
    \end{tikzpicture}

    \caption{Three situations for the location of $v_{i_X}$ and $v_{j_X}$ with respect to the core $Y$.}% of normal tops in~\cref{clm:ell_ny}}
    \label{fig:3sit}
\end{figure}
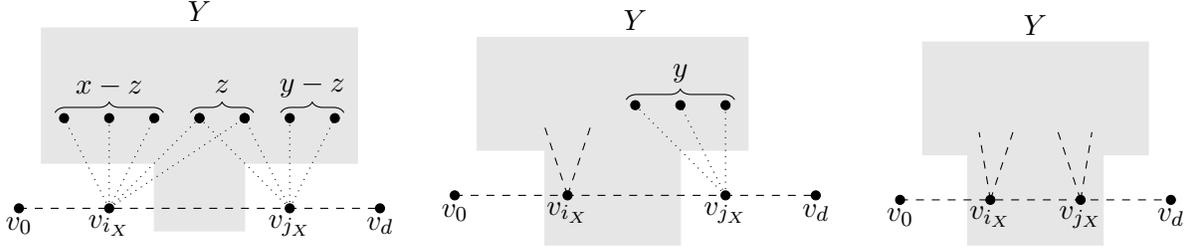

    Now the remainder of the proof is analogous to the final steps of the proof in~\cite{John22}, so we conclude with a summary of those steps.
    The path of a high top $\tau=(T,c)\in \mathcal H$ contains exactly one tree $X$ at distance $1$ from $P$ (where possibly $X=T$).
    We refer to this tree as the ``extension'' of the top $(T,c)$ and denote it by $\ext((T,c))$. Note that every tree $X$ satisfying $d(X,V(P))=1$ is the extension of exactly two high tops, since it lies on one path of each colour.
    Now $\ell((T,\text{red}))\le n-d-1+i_X$ and $\ell((T,\text{blue}))\le n-d-1+(d-j_X)=n-1-j_X$ since one extends the top by adding one vertex at a time starting from a tree with at least one vertex, and $d-i_X$ resp.~$j_X$ vertices are not appearing.
    
    For such a tree $X$ with $d(X,V(P))=1$, there are at least two normal tops $\tau$ for which $\res(\tau)=X$
    (this is verified by considering $i_X=j_X$ and $i_X<j_X$ separately, using $d \ge 2$).
    Let $\nx$ be the collection of all such normal tops $\tau$ for which $\res(\tau)=X$.

    Now we obtain that $\ell\left( \nx \cup\{\tau\in\mathcal H:\ext(\tau)=X\} \right) \le \frac{n-1}{2}$ by combining~\cref{clm:ell_ny} with the inequalities 
    $\frac{d+1}{2} \le \frac{n-1}{2}$, $\abs{\nx} \ge 2$, as well as $\ell((T,\text{red}))\le n-d-1+i_X$ and $\ell((T,\text{blue}))\le n-1-j_X$ for the high top $T$ with $\ext((T,c))=X$. Finally,
    $$2\frac {d+ \mathbf{1}_{j_X>i_X}}2 + n-d-1+i_X + n-1-j_X \le 2(n-1).$$
    Together with $\ell(\mathcal L) =\frac d2 < \frac{n-1}{2}$ and $\ell(\nx) \le \frac{d+1}{2} \le \frac{n-1}{2}$ for forests $X$ that are not trees, this implies that $\ell(\mathcal N \cup \mathcal L \cup \mathcal H) < \frac{n-1}{2}.$

    The average number of vertices in a path associated with a top is thus strictly less than $\frac{n+1}{2}$. So this is also true for the paths of at least one of the colours, without loss of generality red.
    The latter implies that the proportion of subtrees of $G$ containing $v_0$ is greater than $\frac{2}{n+1}.$    
\end{proof}

The proof of~\cref{conj:main}~\ref{itm:min_Pn} given~\cref{prop:localsubtreeorder} and~\cref{Fraction} is now immediate and completely analogous to the proofs in~\cite{John22,Vince22} (see~\cite[Proof of Theorem 1]{John22}, end of page 6, with $A$ replaced by $\mu$).
In essence, by induction and the inequality given by~\cref{Fraction} for a vertex $v$, we have
\begin{align*}
    \mu(G) &=\frac{(N(G)-N(G,v))\mu(G \setminus v) + N(G,v)\mu(G,v)}{N(G)} \\
    &\ge \frac{n+1}{3}+\frac{N(G,v)}{N(G)}\left(\frac{n+1}2-\frac{n+1}{3}\right) \ge\frac{n+1}{3}+\frac{2}{n+1}\cdot\frac{n+1}{6}=\frac{n+2}{3},
\end{align*}
and equality holds if and only if $G$ is a path (which is necessary by~\cref{Fraction}).

\section{The clique maximises the mean subtree order among connected graphs}
\label{sec:maximum_Kn}

In this section, we prove a result on the number of spanning trees and ``almost spanning'' trees that implies that the complete graph has the greatest mean subtree order. Let $s_k(G)$ be the number of subtrees of $G$ of order $k$ as before. In particular, if $G$ has $n$ vertices, then $s_n(G)$ is the number of spanning trees.
    
The mean subtree order can be written as
$$\mu(G) = \frac{\sum_{k=1}^n k s_k(G)}{\sum_{k=1}^n s_k(G)}.$$
Moreover, 
$$p(G)=\frac{ s_n(G)}{N(G)}=\frac{ s_n(G)}{ \sum_{k=1}^n s_k(G)}$$ is the fraction of subtrees of $G$ that are spanning, or equivalently the probability that a random subtree of $G$ is spanning. 
We first prove 

\begin{prop}\label{conj1}
For every graph $G$ with $n\ge 2$ vertices, we have\footnote{Note that $s_n(G) = 0$ if $G$ is disconnected. For this reason, the inequality is not presented in the more symmetric form $\frac{s_{n-1}(G)}{s_n(G)} \ge \frac{s_{n-1}(K_n)}{s_n(K_n)}$.}
$$s_{n-1}(G) \ge \frac{s_{n-1}(K_n)}{s_n(K_n)} s_n(G).$$
\end{prop}

\begin{proof}
For $m\geq1$, let $R_m$ denote the rooted tree shape obtained by
choosing a labelled tree on $[m]$ uniformly at random, rooting it at a uniformly random vertex,
and then forgetting the labels.  By
\cite[Theorem~2.1]{LyonsPeledSchramm} (which builds on~\cite{LuczakWinkler}), $R_{n-1}$ and $R_n$ can be
coupled so that $R_{n-1}$ is a rooted subtree of $R_n$ with the same root.  Since their
orders differ by $1$, $R_n$ is obtained from $R_{n-1}$ by adjoining a
single leaf $u$.

Fix a set $W$ of size $n$, and label the vertices of $R_n$ uniformly with the elements of $W$.  Let
$w$ be the label of $u$, and let $T$ be the resulting unrooted labelled
tree.  Then $T$ is uniform on the set of trees whose vertices are labelled with the elements of $W$.  Moreover, $w$ is uniform on
$W$, and, conditional on $w$, the inherited labeling of $R_{n-1}$ is
uniform on $W\setminus\{w\}$.  Since $R_{n-1}$ has the shape of a
uniformly rooted uniform labelled tree, after forgetting its root we
conclude that $T-u$ is uniform among the trees with labels
$W\setminus\{w\}$.  

Now take $W$ to be the vertex set of $G$ (more precisely, label the vertices of $G$ with the elements of $W$, and identify each vertex with its label). Since $T$ is a uniformly random spanning tree on the complete graph with vertex set $W$, the probability for $T$ to be a subgraph (thus a spanning tree) of $G$ is
$$\frac{s_n(G)}{s_n(K_n)}.$$
Similarly, since $T-u$ is a uniformly random subtree of order $n-1$ of the complete graph with vertex set $W$, the probability for $T - u$ to be a subgraph of $G$ is
$$\frac{s_{n-1}(G)}{s_{n-1}(K_n)}.$$
If $T$ is a subgraph of $G$, then $T-u$ is automatically also a subgraph of $G$. So it follows that
$$\frac{s_{n-1}(G)}{s_{n-1}(K_n)} \geq \frac{s_n(G)}{s_n(K_n)},$$
which readily implies~\cref{conj1}.
\end{proof}

Now, we derive the main results as corollaries of~\cref{conj1}.

\begin{lem}\label{lem:generalratio}
For every graph $G$ with $n$ vertices and every positive integer $k < n$, we have
$$s_k(G) \ge \frac{s_k(K_n)}{s_{k+1}(K_n)} s_{k+1}(G).$$
\end{lem}

\begin{proof}
We have
$$s_k(G) = \sum_{\substack{A \subseteq V(G) \\ \abs{A} = k}} s_k(G[A]) = \frac{1}{n-k} \sum_{\substack{B \subseteq V(G) \\ \abs{B} = k+1}} \sum_{v \in B} s_k(G[B-v]) = \frac{1}{n-k} \sum_{\substack{B \subseteq V(G) \\ \abs{B} = k+1}} s_k(G[B]).$$
To see why the second identity holds, note that every subset $A$ of $k$ vertices is contained in precisely $n-k$ subsets $B$ of $k+1$ vertices and thus counted $n-k$ times in the double sum. Now we can apply \cref{conj1} to every induced subgraph of the form $G[B]$:
$$s_k(G) \ge \frac{1}{n-k} \sum_{\substack{B \subseteq V(G) \\ \abs{B} = k+1}} \frac{s_{k}(K_{k+1})}{s_{k+1}(K_{k+1})} s_{k+1}(G[B]) = \frac{1}{n-k} \cdot \frac{s_{k}(K_{k+1})}{s_{k+1}(K_{k+1})} s_{k+1}(G).$$
Using the fact that every induced subgraph of a complete graph is again complete, one can apply the same reasoning as before in the opposite direction to show that
$$\frac{1}{n-k} \cdot \frac{s_{k}(K_{k+1})}{s_{k+1}(K_{k+1})} = \frac{s_k(K_n)}{s_{k+1}(K_n)}.$$
Alternatively, one can also use the explicit formula $s_k(K_n) = \binom{n}{k} k^{k-2}$. The lemma follows immediately.
\end{proof}

\begin{lem}\label{lemma:quotients}
Let $j,k,n$ be positive integers with $j \le k \le n$. For every graph $G$ with $n$ vertices, we have
$$s_j(G) \ge \frac{s_j(K_n)}{s_k(K_n)} s_k(G).$$
\end{lem}

\begin{proof}
From~\cref{lem:generalratio}, we obtain
$$
s_j(G)\ge \frac{s_j(K_n)}{s_{j+1}(K_n)} s_{j+1}(G) \ge \frac{s_{j}(K_n)}{s_{j+2}(K_n)} s_{j+2}(G) \ge \cdots \ge \frac{s_{j}(K_n)}{s_{k}(K_n)} s_{k}(G).
$$
\end{proof}

In~\cite{CGMV18}, the authors list the question of determining the maximum of $p(G)$ over the set of all graphs with a given number of vertices and edges as a problem for further study, and even when only the number of vertices is prescribed, this was an open problem. The latter problem is resolved below; the problem with both order and size prescribed remains open.

\begin{cor}\label{conj2}
For every graph $G$ with $n$ vertices, we have $p(G) \le p(K_n)$.
\end{cor}

\begin{proof}
By~\cref{lemma:quotients}, we have
$$s_j(G) \ge \frac{s_j(K_n)}{s_n(K_n)} s_n(G)$$
for all $j \in \{1,2,\ldots,n\}$. Summing over all $j$, we obtain the following, which is equivalent to $p(G) \le p(K_n)$:
$$s_1(G) + s_2(G) + \cdots + s_n(G) \ge \frac{s_1(K_n)+s_2(K_n) + \cdots + s_n(K_n)}{s_n(K_n)} s_n(G).$$
\end{proof}

\begin{remark}
 As a side note to~\cref{conj2}, it is easy to show that $p(G)$ is minimised (among connected graphs of order $n$) by the star $S_n$.
 More precisely, $p(G) \ge p(S_n)$ for every connected graph $G$ of order $n$, with equality if and only if $G$ is a star.
 To see this, note that the star is the tree with the greatest number of subtrees, namely $2^{n-1}+n-1$ (see \cite[Theorem 3.1]{SW05}). This implies that every spanning tree of $G$ has at most $2^{n-1}+n-1$ subtrees. Every subtree of $G$ is contained in a spanning tree, thus the total number of subtrees of $G$ is less than or equal to $(2^{n-1}+n-1)s_n(G)$. Equality can only hold if every spanning tree is a star and every subtree belongs to only one spanning tree (in particular for a singleton subtree this implies that $G$ itself is a tree), which is only the case if $G$ itself is a star. The statement follows.
\end{remark}

Finally, we prove~\cref{conj:main}~\ref{itm:max_Kn}, which we restate for the convenience of the reader.

\begin{cor}\label{conj3}
For every graph $G$ with $n$ vertices, we have $\mu(G) \le \mu(K_n)$, with equality only attained by $K_n.$
\end{cor}

\begin{proof}
Fix $j \in \{1,2,\ldots,n-1\}$. \Cref{lemma:quotients} gives us
$$s_a(K_n) s_b(G) \le s_a(G) s_b(K_n)$$
for all positive integers $a$ and $b$ with $1 \le a \le j < b \le n$. Summing these inequalities over all $a$ and $b$, we obtain
$$\sum_{a=1}^j s_a(K_n) \sum_{b=j+1}^n s_b(G) \le \sum_{a=1}^j s_a(G) \sum_{b=j+1}^n s_b(K_n).$$
Add the double sum $\sum_{a=j+1}^n s_a(G) \sum_{b=j+1}^n s_b(K_n)$ on both sides:
$$\sum_{a=1}^n s_a(K_n) \sum_{b=j+1}^n s_b(G) \le \sum_{a=1}^n s_a(G) \sum_{b=j+1}^n s_b(K_n).$$
This inequality trivially holds (with equality) for $j=0$ as well. Summing over all $j$ (from $0$ to $n-1$) gives us
$$\sum_{a=1}^n s_a(K_n) \sum_{b=1}^n b s_b(G) \le \sum_{a=1}^n s_a(G) \sum_{b=1}^n b s_b(K_n),$$
which is equivalent to the inequality in \cref{conj3}.
Since for $a=1$ and $b=2$, whenever $G \not= K_n$, $s_1(K_n)s_2(G)< s_1(G)s_2(K_n),$ the inequality $\mu(G) \le \mu(K_n)$ is strict whenever $G\not=K_n.$
\end{proof}

\section{Counterexamples to the universal version of~\cref{conj:main2}}\label{sec:nonuniversal}

In this section, we consider the change in the mean subtree order under adding or removing an edge. Specifically, in the first part we construct examples of graphs with the property that ``many'' potential edge additions reduce the mean subtree order, while in the second part we construct examples of graphs with the property that ``many'' potential edge removals increase the mean subtree order.

\subsection{Proof of~\cref{thr:strong_counterexamples_toConj7.4Chin}}

Here we provide a construction of graphs $G$ for which there is essentially only one edge satisfying $\mu(G+e)>\mu(G).$
Let us first collect a few facts about subtrees of complete graphs that will be useful in the following.

\begin{itemize}
\item The number of spanning trees of a complete graph $K_n$ is well known to be $n^{n-2}$, and it was shown that the number of subtrees asymptotically only differs by a constant factor: $N(K_n) \sim e^{1/e} n^{n-2}$ (see~\cite[Thm. 3.3]{CGMV18} as well as \cite{CGMV15}).
\item By~\cite[Cor.~2]{Wagner21_prob_spanning}, $\mu(K_n) = n - \frac1e + o(1)$.
\item The number of subtrees that contain a specific vertex $v$ is
$$N(K_n) - N(K_{n-1}) \sim N(K_n) \sim e^{1/e} n^{n-2}.$$
\item Let us fix a second vertex $w$. Among the subtrees that contain $v$, precisely $N(K_{n-1}) - N(K_{n-2}) \sim N(K_{n-1}) \sim e^{1/e} (n-1)^{n-3} \sim e^{1/e-1} n^{n-3}$ do not contain $w$. It follows that
$$\mu(K_n,v) = n - (n-1) \frac{N(K_{n-1}) - N(K_{n-2})}{N(K_n) - N(K_{n-1})} = n - \frac1e + o(1).$$
\item Every edge $uv$ is contained in the same number of subtrees. Since the average number of edges in a subtree is $\mu(K_n) - 1 = n - O(1)$, it follows that the proportion of subtrees that contain $uv$ is
$$\frac{n-O(1)}{\binom{n}{2}} = \frac{2}{n} + O\big(n^{-2}\big),$$
which means that the number of subtrees containing a fixed edge $uv$ is asymptotically equal to $\frac{2}{n} N(K_n) \sim 2e^{1/e} n^{n-3}$.
\item If we fix a third vertex $w \notin \{u,v\}$, there are (by the same argument applied to $K_{n-1}$) asymptotically $\frac{2}{n-1} N(K_{n-1}) \sim 2e^{1/e-1} n^{n-4}$ subtrees that contain $uv$, but not $w$. Hence we also have
$$n - \mu(K_n,uv) \sim (n-2) \cdot \frac{\frac{2}{n-1} N(K_{n-1})}{\frac{2}{n} N(K_n)} \sim \frac{1}{e}$$
or
$$\mu(K_n,uv) = n - \frac{1}{e} + o(1).$$
Thus $\mu(K_n)$, $\mu(K_n,v)$ and $\mu(K_n,uv)$ are equal up to $o(1)$, but as~\cref{prop:VET_mue>muv>mu} in the following section shows (more generally for graphs that are both vertex- and edge-transitive), we have $\mu(K_n) < \mu(K_n,v) < \mu(K_n,uv)$ for all $n \ge 2$.
\end{itemize}

Now let $D_{n,w}$ be a barbell of order $n$, which consists of two cliques of order $w$ and a path that connects the two cliques (\cref{fig:D_{14,6}graph} depicts a concrete example, namely $D_{14,6}$). In the following, we show that barbells provide the desired examples for a suitable choice of $n$ and $w$.

\begin{prop}
 % For fixed $k$, 
 Let $G=D_{n,w}$ with $w=o(n)$ and $w=\omega( \log n)$ be a barbell of sufficiently large order $n$.
 % Then one can add $k$ edges to $G$, such that the mean subtree order decreases with roughly $\frac n3.$
 There is essentially only one edge ($2w-2$ distinct choices, all belonging to the same isomorphism class) whose addition increases the mean subtree order.
\end{prop}

\begin{proof}
 % Form $G_k$, with $\mu(G_k)\sim \frac 23n$ and 
 
 Let the internal vertices (vertices of degree $2$) on the path be $v_1, v_2, \ldots, v_{n-2w}.$ Let $v_0$ and $v_{n-2w+1}$ be the ends of the path, i.e., the two vertices of maximum degree in $G$.
 Let $u_1, u_2, \ldots, u_{w-1}$ and $u'_1, u'_2, \ldots, u'_{w-1}$ be the other vertices of the two cliques, respectively. 
% Note that almost all subtrees of $G$ contain the whole path and $\mu(G) \sim n-\frac{2}{e}$ (by~\cite[Cor.~2]{Wagner21_prob_spanning}).

For our choice of $w$, $w^{w-2}$ grows faster than $n$ (in fact faster than any power of $n$), which implies that almost all subtrees of $G$ contain the entire path $v_0v_1 \cdots v_{n-2w+1}$. All other subtrees are negligible. It follows from our observations on subtrees of complete graphs that $\mu(G) = n - \frac{2}{e} + o(1)$.

 To prove the inequality $\mu(G+e)<\mu(G),$ it is sufficient to show that $\mu(G+e,e)<\mu(G)$, since $\mu(G+e)$ is a weighted average of $\mu(G+e,e)$ and $\mu(G)$---note here that $\mu(G)$ can be interpreted as the average subtree order of all subtrees in $G+e$ that do not contain $e$.
 We do so for all possibilities (up to relabelling, i.e., for a representative of each isomorphism class). Here the cases are ordered by the proportion of edges that belong to each case.
 
 \textbf{Case 1:} $e=v_i v_j$ with $\abs{i-j}\ge 2.$
 
It is still the case (by the same reasoning) that all but a negligible proportion of the subtrees of $G+e$ that contain $e$ also contain $v_0$ and $v_{n-2w+1}$.  The subtrees containing the edge $e=v_i v_j$ are missing $\frac{\abs{i-j}-1}3$ of the vertices in $\{v_i, v_{i+1}, \ldots, v_{j}\}$ on average. To see this, observe that the complement of a subtree that contains $e$ necessarily induces a subpath between $v_i$ and $v_j$.
 This excluded subpath (which cannot be a single vertex) has an average length of $\frac{\abs{i-j}+2}3$ and thus an average number of internal vertices equal to $\frac{\abs{i-j}-1}3$.
 % and between $v_i$ and $v_j$ there is a loss of $\frac{\abs{i-j}-1}3$ on average (a subpath between $v_i$ and $v_j$ is excluded) if the subtree contains $e$ compared with containing the path $v_iv_{i+1} \ldots v_j.$
 This implies that $\mu(G+e,e)=\mu(G)-\frac{\abs{i-j}-1}3+o(1)$. If $n$ is sufficiently large, we can conclude that $\mu(G+e,e)<\mu(G).$
 % This difference is substantial to conclude.
 
 \textbf{Case 2:} $e=u_1 v_i$ with $i\ge 2.$ 

In this case, we can use the same argument as before to show that all but a negligible proportion of the subtrees of $G+e$ that contain $e$ have to contain the entire path $v_i v_{i+1} \cdots v_{n-2w+1}$ as well, and contain on average all but $\frac1e + o(1)$ of the vertices of the clique induced by the $u_i'$s and $v_{n-2w+1}$. For the remaining part, we have the following three possibilities:
\begin{itemize}
\item The subtree contains $v_0$, but not the entire path from $v_0$ to $v_i$. By the inclusion-exclusion principle, there are $N(K_w) - 2N(K_{w-1}) + N(K_{w-2}) \sim N(K_w)$ possibilities to choose the subtree induced by the clique that consists of the $u_i$s and $v_0$, and $\binom{i+1}{2}$ choices for the forest induced on $v_0v_1 \cdots v_i$. The former part misses $\frac1e + o(1)$ of the vertices on average, the latter $\frac{i-1}{3}$.
\item The subtree contains the entire path from $v_0$ to $v_i$. The rest is a forest with two components inside the clique, where one component contains $v_0$, the other $u_1$. By adding the edge $v_0u_1$, we obtain a subtree of the clique that contains this edge, and we know that these subtrees miss $\frac{1}{e} + o(1)$ vertices on average.
\item The subtree does not contain $v_0$. The rest consists of a subtree of a clique of order $w-1$ that contains the fixed vertex $u_1$ and a path of the form $v_jv_{j+1} \cdots v_i$ for some $j \in \{1,2,\ldots,i\}$. In total, these subtrees miss $\frac{i+1}{2} + \frac1e + o(1)$ vertices on average.
\end{itemize}
Now note that the second and third case are asymptotically negligible compared to the first: there are $\sim \binom{i+1}{2} \cdot e^{1/e} w^{w-2}$ subtrees in the first case, $\sim 2e^{1/e} w^{w-3}$ in the second, and $\sim i e^{1/e-1} w^{w-3}$ in the third.

Putting everything together, we have $\mu(G+e,e)=n - \frac{i-1}{3} - \frac{2}{e} + o(1) = \mu(G)-\frac{i-1}3+o(1)<\mu(G)$ in this case if $n$ is sufficiently large.
 
 \textbf{Case 3:} $e=u_1 u'_1.$
 
 In this case, almost all subtrees containing $e$ contain $v_0$ and $v_{n-2w+1}$, but not the whole path $v_0 \ldots v_{n-2w+1}$ (the precise verification is similar to the previous two cases). It follows that $\mu(G+e,e) \sim \frac{2n}3$, which is smaller than $\mu(G)$.

 \textbf{Case 4:} $e=u_1 v_1.$

 This is the only case where $\mu(G+e,e)>\mu(G)$ and thus $\mu(G+e) > \mu(G)$. To see this, we distinguish different types of subtrees of $G$ and $G+e$:
 \begin{itemize}
    \item Type A: subtrees of $G$ that do not contain the edge $v_0v_1$,
    \item Type B: subtrees of $G$ that contain the edge $v_0v_1$,
    \item Type C: subtrees of $G+e$ that contain $e = u_1v_1$, but not $v_0v_1$,
    \item Type D: subtrees of $G+e$ that contain both $e = u_1v_1$ and $v_0v_1$.
 \end{itemize}
Note that $\mu(G)$ is obtained as an average over subtrees of types A and B, while $\mu(G+e,e)$ is obtained as an average over subtrees of types C and D. There is a straightforward bijection between types B and C (replacing $v_0$ by $u_1$ and vice versa), so those types have the same average order. It has already been established that $\mu(G) = n - \frac{2}{e} + o(1)$, while the average order of subtrees of type A is clearly less than $n-w$, thus less than the average order of subtrees of type B (or C). So if we can show that subtrees of type D have greater average order than those of type B (or C), we are done with the proof that $\mu(G+e,e)>\mu(G)$.

Every subtree of type B is obtained by independently choosing two parts and combining them:
\begin{itemize}
    \item a subtree of the clique induced by the $u_i$s and $v_0$ that contains $v_0$,
    \item a subtree of the graph induced by the $u_i'$s and $v_1,\ldots,v_{n-2w+1}$ that contains $v_1$.
\end{itemize}
The same is true for subtrees of type D, except that for the first part, a two-component forest of the clique is chosen whose components contain $v_0$ and $u_1$ respectively. Equivalently, one can choose a subtree of the clique that contains the edge $v_0u_1$ and remove that edge. It follows from~\cref{prop:VET_mue>muv>mu} in the following section that subtrees of a complete graph that contain a fixed edge have more vertices on average than subtrees containing a fixed vertex. Hence we can conclude that the average order of subtrees of type D is indeed greater than that of subtrees of type B, completing the proof. 
\end{proof}

\begin{remark}
The increase in the mean subtree order in Case 4 is actually tiny: a more precise asymptotic analysis shows that $\mu(G + u_1v_1) - \mu(G) \sim \frac{1}{ew^3}$.
\end{remark}

The following conjecture was communicated privately by Yanli Hao: ``For every maximal matching $M$ of the complement graph of $G$, adding $M$ to $G$ increases the average order of subtrees.'' We disprove this in a strong sense with the following counterexample.

\begin{ex}
 The graph $G=D_{14,6}$ has the property that for every maximal matching $M$ in its complement, $\mu(G+M)<\mu(G)$. Up to isomorphism, there is only one edge $e \in G^c$ for which $\mu(G+e) >\mu(G)$.
\end{ex}

%\begin{proof}
 This has been verified by computer by computing the mean subtree order for all choices of $G+e$ and $G+M$ (taking isomorphisms into account). The source code can be downloaded from \url{https://github.com/JorikJooken/averageSubtreeOrderOfGraphs}.
%\end{proof}

\begin{figure}[ht]
\centering
\begin{tikzpicture}[scale=1.1]

\node at (1.65,0.3) {$v_0$};
\node at (5.85,0.3) {$v_3$};
\node at (3.25,-0.25) {$v_1$};
\node at (4.75,-0.25) {$v_2$};
\node at (0.75,1.65) {$u_1$};
\node at (6.75,1.65) {$u'_1$};

\draw[thick, dashed, red] (60:1.5)--(6,0);
\draw[thick, dashed, red] (60:1.5)--(4.5,0);
\draw[thick, dashed, red] (60:1.5)--(3,0);
\draw[thick, dashed, red] (60:1.5)--(6.75,{1.5*sin(60)});

\draw[thick, dashed, red] (1.5,0) arc (240:300:3cm);
\draw[thick, dashed, red] (1.5,0) arc (240:300:4.5cm);

\foreach \x in {0,60,...,300}{\draw[fill] (\x:1.5) circle (0.15);
\foreach \a in {60,120,180}{
\draw(\x+\a:1.5) -- (\x:1.5);
}
}

\foreach \x in {0,60,...,300}{
\draw[fill] ({7.5+1.5*cos(\x)},{1.5*sin(\x)}) circle (0.15);
\foreach \a in {60,120,180}{
\draw ({7.5+1.5*cos(\x)},{1.5*sin(\x)}) -- ({7.5+1.5*cos(\x+\a)},{1.5*sin(\x+\a)});
}
}
\draw (1.5,0)--(6,0);

\draw[fill] (3,0) circle (0.15);
\draw[fill] (4.5,0) circle (0.15);

\end{tikzpicture}
\caption{The graph $D_{14,6}$ and $6$ possible edges to add.}\label{fig:D_{14,6}graph}
\end{figure}
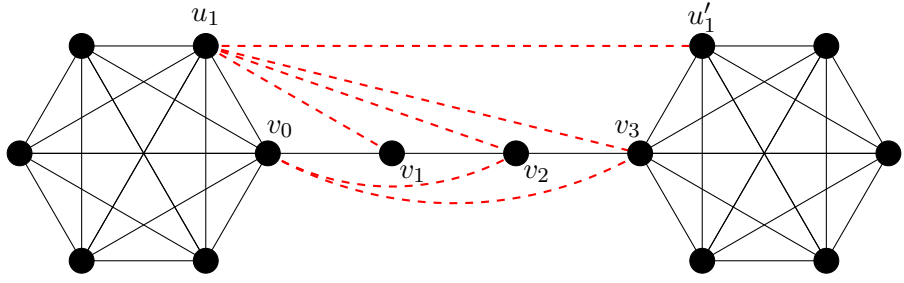

\subsection{Deleting an edge can increase the mean subtree order}\label{subsec:deletingedge}

As we have seen that adding an edge can decrease the mean subtree order, deleting an edge can conversely increase it. Here we illustrate this with concrete examples of graphs that have many edges for which this is the case.

The join of two graphs $G$ and $H$, written as $G \vee H$, is defined as the graph with vertex set $V(G) \cup V(H)$ and edge set $E(G) \cup E(H) \cup \{ xy \mid x \in V(G), y \in V(H) \}$. Let $G=K_n \vee m K_1$ (equivalently, $K_{n+m}$ with a $K_m$ removed).
This is a complete multipartite graph with partition classes of size $m$ and $n$ times $1$.
Note that every induced subgraph $G'$ of $G$ is itself of the form $K_{n'} \vee m' K_1.$

% With Kirchhoff's matrix tree theorem, we can compute the number of spanning trees in the graph $G$ (and thus also every $G'$).
% For this, we compute the determinant of the Laplacian matrix $L$ of $G$, with one row and one column removed.
% $$\tilde{L}(G)=\begin{bmatrix}
%  (n+m)I_{n-1}-J_{n-1} & -J_{n-1,m} \\
%  -J_{m,n-1} & n I_m \\
% \end{bmatrix}$$

% Here $J_{k,d}$ denotes an $k \times d$-matrix for which all its entries equal $1$. With $I_n$ and $J_n$ we denote the $n\times n$ identity matrix and $J_{n,n}$ respectively. So in particular, one submatrix of $L$ equals
% $$(n+m)I_{n-1}-J_{n-1}=
% \begin{bmatrix}
%  n+m-1 & -1 & \cdots & -1 & -1 \\
%  -1 & n+m-1 & \cdots & -1 & -1 \\
%  \vdots & \vdots & \vdots & \ddots & \vdots \\
%  -1 & -1 & \cdots & n+m-1 & -1 \\
%  -1 & -1 & \cdots & -1 & n+m-1 \\
% \end{bmatrix}.$$

% Let $e$ be an edge between two vertices of degree $n+m-1$ in $G$ (an edge of the central $K_n$ within the sunflower $G$).
% Then the Laplacian matrix $L(G \setminus e)$ equals $L(G),$ except on four entries.
% We can take $\tilde L(G \setminus e)$ equal to $ \tilde L(G),$ except that its leftupper element $\tilde L_{1,1}$ equals $n+m-2$ instead of $n+m-1.$

% By row operations, we find that the determinant of $\tilde{L}(G)$ equals $(n+m)^{n-2}n^{m-1}\cdot(n+m)=(n+m)^{n-1}n^{m-1}.$

Let $e$ be an edge between two vertices of degree $n+m-1$ in $G$.
Now $G \setminus e$ is a complete multipartite graph with partition classes of size $m$, $2$ and $n-2$ times $1$.

The number of spanning subtrees of $G$ and $G \setminus e$ is given by $f(n,m)=(n+m)^{n-1}n^{m-1}$ and
$g(n,m)=(n+m-2)(n+m)^{n-2}n^{m-1}$, respectively. These formulas can be derived from Kirchhoff's matrix-tree theorem, see~\cite{Lewis99}. Note in particular that $f(n,0) = n^{n-2}$ is precisely Cayley's formula.

% We will denote this expression as $f(n,m)=(n+m)^{n-1}n^{m-1}.$

% By the difference of the single entry, $\det \tilde{L}(G \setminus e)$ equals $\det \tilde{L}(G)$ minus the determinant of the minor of the leftupper element,
% $$\det \tilde{L}(G \setminus e)=(n+m)^{n-1}n^{m-1}-2(n+m)^{n-2}n^{m-1}=(n+m-2)(n+m)^{n-2}n^{m-1}.$$
% We will use $g(n,m)=(n+m-2)(n+m)^{n-2}n^{m-1}$ for this number.

Now the number of subtrees in $G=K_n \vee m K_1$ can be computed by taking the number of spanning trees over all subgraphs of the form $K_{j} \vee i K_1.$
There are $\binom{n}{j}\binom{m}{i}$ possibilities for this subgraph.
Here, $i,j$ can be arbitrary nonnegative integers, except for the case $j=0$ where the only subtrees are single vertices.
So the number of subtrees of $G$ is
$$N(G)=m+\sum_{j=1}^n \sum_{i=0}^m \binom{n}{j}\binom{m}{i}f(j,i),$$
%Note here that $f(j,0)=j^{j-2}$ (Cayley's formula for the number of spanning trees of a complete graph).
and the total order of all subtrees of $G$ is 
$$R(G)=m+\sum_{j=1}^n \sum_{i=0}^m \binom{n}{j}\binom{m}{i}(i+j)f(j,i). $$

The induced subgraphs of $G \setminus e$ are either of the form $K_{j} \vee i K_1$ or $(K_{j} \vee i K_1)\setminus e.$
There are $\binom{n-2}{j-2}\binom{m}{i}$ possibilities for the latter.
In the same way as for $G$, we find that 
\begin{align*}
    N(G \setminus e)&=m+\sum_{i=0}^m \binom{m}{i}\sum_{j=1}^n  \left( \left(\binom nj-\binom{n-2}{j-2}\right)f(j,i) + \binom{n-2}{j-2}g(j,i) \right)\mbox{ and } \\
    R(G \setminus e)&=m+\sum_{i=0}^m \binom{m}{i}\sum_{j=1}^n (i+j) \left( \left(\binom nj-\binom{n-2}{j-2}\right)f(j,i) + \binom{n-2}{j-2}g(j,i)\right). 
\end{align*}

For many values of $n$ and $m$, we have $\mu(G)=\frac{R(G)}{N(G)}<\frac{R(G\setminus e)}{N(G\setminus e)}=\mu(G \setminus e)$. Cameron and Mol~\cite{CM21} present the special case $n = 2$ and $m \ge 6$. We verified by computer (see~\url{https://github.com/JorikJooken/averageSubtreeOrderOfGraphs}, document DeletionEdge\_mu\_difference) that this also holds for all $n \in [10,400]$ that are multiples of $10$, as well as for $n = 800$, if $m= \frac{7n}{10}+2$. For these examples with $m = \frac{7n}{10} + 2$, roughly $\frac{5}{12}$ of the edges have the property that $\mu(G)<\mu(G \setminus e)$.

% There is a bound $f(n)$, such that for $m\ge f(n),$ $\mu(G)$ increases when deleting an edge $e$ of the clique; 
% Hereby $f(2)=6$ (done by Cameron-Mol) and $f(n)=0.7n+2$ when $10\le n \le 400$ is a multiple of $10$, or $n =800$ (we checked the latter claim in.
% This implies that for these examples, roughly $\frac{5}{12}$ of the edges have the property that $\mu(G)<\mu(G \setminus e).$

%For order up to $9$, these constructions looked like the ones with the largest proportion of edges whose removal increases the mean subtree order.

\section{The mean subtree order is not monotone in general graphs}\label{sec:nomonotonicity}

In~\cite{jamison_monotonicity_1984}, it was proven that the average order of the subtrees of a tree $T$ that contain a particular subtree~$S$, denoted by $\mu(T, S)$, is increasing in $S$: $\mu(T,S) < \mu(T,S')$ if $S$ is (properly) contained in $S'$. In particular, if $e$ is an edge incident with a vertex $v$, we have $\mu(T,e) > \mu(T,v)$.
Moreover, the local mean $\mu(T,v)$ at a vertex $v$ is always greater than or equal to the global mean $\mu(T)$, see \cite{jamison_average_1983}.
Nevertheless, as we will see in this section, the following analogues for general graphs do not hold.

\begin{enumerate}
 \item\label{itm:1} For every vertex $v$ in a graph $G$, $\mu(G,v) \ge \mu(G).$
 \item\label{itm:2} For subtrees $S \subset S'\subset G,$
 $\mu(G,S)<\mu(G,S')$.
 \item\label{itm:3} For every vertex $v$ in a graph $G$, there is an edge $e$ incident with $v$ for which $\mu(G,e)>\mu(G,v).$
\end{enumerate}

If~\ref{itm:1} (or a weaker version guaranteeing existence of a non-cut vertex $v$ such that $\mu(G,v) \ge \mu(G)$) and~\ref{itm:3} were true, there would be an edge $e$ for which $\mu(G,e)>\mu(G)$ and thus $\mu(G \setminus e)<\mu(G),$ implying~\cref{conj:main2}~\ref{itm:main2_Pn}.
However, as the three statements above are all false, one cannot argue in this way. We will show this by constructing explicit counterexamples.

Let $D^{\star a}_{n,w}$ be a modified barbell, which is constructed by taking a barbell $D_{n-a,w}$ and connecting the two vertices of maximum degree by a path of length $a+1$.
Similarly, let $DB_{n,w}$ be a double broom of order $n$, which consists of two stars of order $w$ whose centres are connected by a path.
The graph $DB^{\star a}_{n,w}$ is obtained from $DB_{n-a,w}$ by connecting the two maximum degree vertices by an additional path of length $a+1$.

\begin{prop}\label{prop:non-monotone}
 There exists a pair $(G,v)$ consisting of a graph $G$ and a vertex $v \in V(G)$ for which $\mu(G)>\mu(G,v)>\mu(G,e)$ holds for every edge $e$ containing $v$.
\end{prop}

\begin{proof}
Consider the graph $G = DB^{\star 1}_{n,w}$, where $w=\omega(\log n)$ and $w=o(n)$, and $n$ is sufficiently large. It can be seen as a cycle $C_{n-2w+2}$ with vertex set $V=\{v_0,v_1,v_2, \ldots, v_{n-2w+1}\}$ (indices in cyclic order), where $v_0$ and $v_2$ are each connected with $w-1$ additional vertices (thus $G \setminus v_1$ is a double broom). The number of subtrees containing neither $v_0$ nor $v_2$ is $O(n^2)$, the number of subtrees containing $v_0$, but not $v_2$ (or vice versa) is $O(n2^w)$. Since there are clearly more than $2^{2w-2}$ subtrees that contain both $v_0$ and $v_2$, and $n = o(2^w)$
by the choice of $w$, it follows that almost all subtrees contain $v_0$ and $v_2$.
 There are $\Theta(n^2)$ subtrees of the cycle which contain the path $v_0v_1v_2$, while there are only three subtrees containing the path $v_2v_3 \ldots v_{n-2w+1}v_0$. Thus almost all subtrees consist of the path $v_0v_1v_2$, a collection of leaves, and paths of the form $v_2v_3 \cdots v_p$ and $v_qv_{q+1} \cdots v_{n-2w+1}v_0$ for some $p < q$. In $G \setminus v_1$, almost all subtrees consist of the path $v_2v_3 \ldots v_{n-2w+1}v_0$ and a collection of leaves.
  From this, one can conclude that $\mu(G)= \frac{2n}{3} + o(n)$ while $\mu(G \setminus v_1) = n-o(n)$, thus $\mu(G) > \mu(G,v_1)$.

 For the second part, we can take $e=v_0v_1$ and $v=v_1$ by symmetry.
 Now $\mu(G,v)>\mu(G,e)$ is equivalent to
 $\mu(G \setminus e,v )>\mu(G,v)$, since $\mu(G,v)$ is a weighted average of $\mu(G \setminus e,v )$ and $ \mu(G,e)$.
 Since $\mu(G \setminus e,v ) = n - o(n)$ and $\mu(G,e)=\frac{2n}{3} + o(n)$, the conclusion is also immediate.
\end{proof}

\begin{figure}[ht]
\centering

\begin{tikzpicture}[scale=0.8]

\node at (1.65,-0.3) {$v_0$};
\node at (5.85,-0.3) {$v_2$};
\node at (3.75,-0.3) {$v_1$};

\foreach \x in {150,160,...,210}{
\draw[fill] ({1.5+2.5*cos(\x)},{2.5*sin(\x)}) circle (0.15);
\draw (1.5,0) -- ({1.5+2.5*cos(\x)},{2.5*sin(\x)});
\draw[fill] ({6-2.5*cos(\x)},{2.5*sin(\x)}) circle (0.15);
\draw (6,0) -- ({6-2.5*cos(\x)},{2.5*sin(\x)});
}

\draw (1.5,0)--(6,0);

\draw[fill] (3.75,0) circle (0.15);
%\draw[fill] (4.5,0) circle (0.15);

 \draw[] (6,0) arc (0:180:2.25) ;

\foreach \x in {0,1,...,7}{
\draw[fill] ({3.75+2.25*cos(\x*180/7)},{2.25*sin(\x*180/7)}) circle (0.15);
}

\end{tikzpicture}\quad
\begin{tikzpicture}[scale=0.8]

\node at (1.65,-0.3) {$v_0$};
\node at (5.85,-0.3) {$v_2$};
\node at (3.75,-0.3) {$v_1$};

\foreach \x in {0,72,...,288}{\draw[fill] (\x:1.5) circle (0.15);
\foreach \a in {72,144}{
\draw(\x+\a:1.5) -- (\x:1.5);
}
}

\foreach \x in {36,108,...,324}{
\draw[fill] ({7.5+1.5*cos(\x)},{1.5*sin(\x)}) circle (0.15);
\foreach \a in {72,144}{
\draw ({7.5+1.5*cos(\x)},{1.5*sin(\x)}) -- ({7.5+1.5*cos(\x+\a)},{1.5*sin(\x+\a)});
}
}
\draw (1.5,0)--(6,0);

\draw[fill] (3.75,0) circle (0.15);
%\draw[fill] (4.5,0) circle (0.15);

 \draw[] (6,0) arc (0:180:2.25) ;

\foreach \x in {0,30,...,180}{
\draw[fill] ({3.75+2.25*cos(\x)},{2.25*sin(\x)}) circle (0.15);
}
\end{tikzpicture}
\caption{The graphs $DB^{\star 1}_{23,8}$ and $D^{\star 1}_{16,5}$ satisfy $\mu(G,v_1)<\mu(G)$.}\label{fig:graph}
\end{figure}
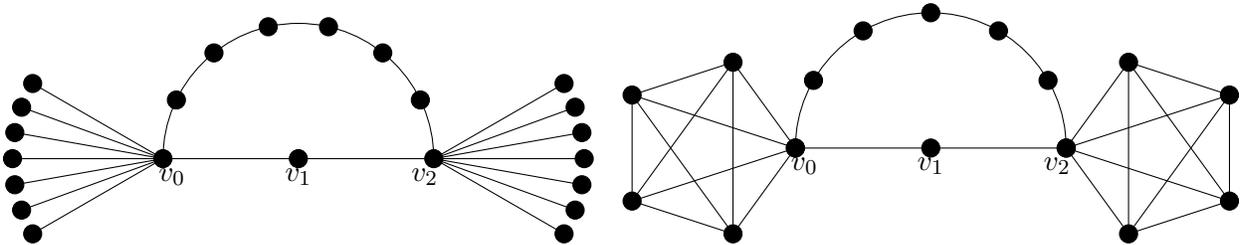

 \begin{remark}
      An alternative would be to take a barbell $D_{n-1,w}$ where $w=\omega(\log n)$ and $w=o(n)$, and to connect its two maximum degree vertices with an additional vertex $v_1$, thus forming $D_{n,w}^{\star 1}$. A similar calculation as for $DB^{\star 1}_{n,w}$ shows that $v_1$ satisfies the statement.
      
      The proposition was also verified explicitly for small concrete examples, see~\url{https://github.com/JorikJooken/averageSubtreeOrderOfGraphs}. Specifically, $DB^{\star 1}_{23,8}$ and $D^{\star 1}_{16,5}$ (a concatenation of $K_5$, $C_8$ and $K_5$) both satisfy the statement. They are presented in~\cref{fig:graph}.
 \end{remark}

\begin{prop}
 For any two fixed positive integers $k$ and $d$, there exist a graph $G$ and subtrees $S \subset S'$ with $\abs S=k$ and $\abs{S'}=k+d$ for which $\mu(G,S)>\mu(G,S')$.
\end{prop}

\begin{proof}
 Take a double broom $DB_{n-d+1,w}$ where $w=\omega(\log n)$ and $w=o(n)$, and $n$ is sufficiently large. Let $S$ be a subpath of $DB_{n-d+1,w}$ of order $k$ containing one of the maximum degree vertices (and no leaves).
 Add a path $P'$ of length $d$ between the two maximum degree vertices  to obtain $G=DB_{n,w}^{\star d-1}.$
 Let $S'$ be the concatenation of $S$ and $P'$. Then we conclude as in the proof of~\cref{prop:non-monotone}: $\mu(G, S') = \frac 23n + o(n)$, while the mean subtree order of the subtrees containing $S$ but not $P'$ is $n-o(n).$
\end{proof}

We remark that for every connected graph $G$ whose order is at least $2$, there is always a vertex $v \in V(G)$ for which the local mean subtree order is larger than the global one, i.e., $\mu(G,v)>\mu(G)$, as a direct application of~\cite[Thm.~3.1]{ARW20} (restated below and with a short proof for completeness). We write $[n]=\{1,2,\ldots,n\}$ for the standard set of cardinality $n$, and $2^{[n]}$ for its power set.
For a family $\A \subset 2^{[n]}$, we write $\A_i=\{A \in \A \mid i \in A\}.$
%https://arxiv.org/pdf/1807.08290.pdf
\begin{lem}
Let $\mathcal A=(A_t)_{t\in I}$ be a non-empty finite indexed
family of subsets of $[n]$, and set
 $\av(\A)=\frac{1}{\abs \A} \sum_{A\in \A} \abs A$.
 If $\A \subset 2^{[n]}$ is a non-uniform family (i.e., not all cardinalities of the sets in $\A$ are the same), then there exists an element $i \in [n]$ such that 
 $\av(\A_i)>\av(\A)$.
\end{lem}

\begin{proof}
We have $$\max \{ \av(\A_i) \mid i \in [n]\} \ge \frac{ \sum_{i=1}^n \sum_{A\in \A_i} \abs A }{\sum_{i=1}^n \abs {\A_i} },$$
since a weighted average (with positive coefficients) is bounded by the maximum.
Writing the cardinalities of the sets $\{(A,i,j)\mid A \in \A; i,j \in A \}$ and $\{(A,i)\mid A \in \A; i \in A \}$ in two different ways, we get
$\sum_{i=1}^n \sum_{A\in \A_i} \abs A = \sum_{A\in \A} \abs A ^2$ and $\sum_{i=1}^n \abs {\A_i} = \sum_{A\in \A} \abs A$. Finally, $\sum_{A\in \A} \abs A ^2 > \frac{1}{\abs \A} \big( \sum_{A\in \A} \abs A \big)^2$ by the Cauchy--Schwarz inequality (strict inequality holds by the assumption that the cardinalities of the sets in $\A$ are not all the same). Putting these together, we obtain
$$\max \{ \av(\A_i) \mid i \in [n]\} \ge \frac{ \sum_{i=1}^n \sum_{A\in \A_i} \abs A }{\sum_{i=1}^n \abs {\A_i} } 
 = \frac{ \sum_{A\in \A} \abs A ^2 }{ \sum_{A\in \A} \abs A} > \frac{\sum_{A\in \A} \abs A }{\abs \A} = \av(\A),$$
 and the statement follows.
\end{proof}
If we apply this lemma to the set $\A$ of all subtrees of $G$, the inequality $\mu(G,v) > \mu(G)$ for some vertex $v$ follows. Similarly, as observed in~\cite{CM21}, there always exists an edge $e$ for which $\mu(G,e)>\mu(G)$.
This implies that an edge-minimal counterexample to~\cref{conj:main}~\ref{itm:min_Pn} must have a bridge.
Now it would be sufficient to prove~\cref{conj:contraction} for $e$ being a bridge\footnote{We remark that~\cite[Eq.~(2.7)]{Ruoyu24} is not true for general graphs (a concrete counterexample is $DB^{\star 1}_{23,8}$), and so a modified strategy is needed.}.

For vertex-transitive graphs, $\mu(G,v)$ is trivially the same for all vertices $v$. Thus the existence of a vertex with $\mu(G,v) > \mu(G)$ implies that this is in fact true for all vertices.

\begin{cor}\label{cor:VT_muv>mu}
    In a vertex-transitive (connected) graph $G$ whose order is at least $2$, $\mu(G,v)> \mu(G)$ for every $v \in V(G)$.
\end{cor}
Along the same lines, one can prove the following analogue for edge-transitive graphs.

\begin{cor}\label{cor:ET_mue>mu}
    In an edge-transitive (connected) graph $G$ (whose order is at least $2$), $\mu(G,e)> \mu(G)$ for every $e \in E(G)$.
\end{cor}

Finally, we consider graphs that are both vertex- and edge-transitive.

\begin{prop}\label{prop:VET_mue>muv>mu}
    In a (connected) graph $G$ (whose order is at least $2$) that is both vertex- and  edge-transitive, $\mu(G,e)> \mu(G,v)>\mu(G)$ for every $e \in E(G)$ and every $v\in V(G)$.
\end{prop}

\begin{proof}
    Recall that $s_k(G)$ denotes the number of subtrees of $G$ with $k$ vertices. 
    Since $G$ is vertex-transitive, every vertex is contained in the same number of subtrees of order $k$ for every $k$, which must be $\frac{k}{n} s_k$ by a simple double-counting argument. Analogously, every edge is contained in $\frac{k-1}{m} s_k$ subtrees of order $k$, where $m=\| G \|$ is the number of edges of $G$. Thus we have
    \begin{align*}
        \mu(G) &= \frac{ \sum_{k=1}^n s_k(G)\cdot k}{ \sum_{k=1}^n s_k(G)} = \frac{ \sum_{k=1}^n k\cdot s_k(G)}{ \sum_{k=1}^n s_k(G)},\\
        \mu(G,v)  &= \frac{ \sum_{k=1}^n k/n\cdot s_k(G)\cdot k}{ \sum_{k=1}^n k/n\cdot s_k(G)} = \frac{ \sum_{k=1}^n k^2\cdot s_k(G)}{ \sum_{k=1}^n k\cdot s_k(G)},\\
        \mu(G,e)&= \frac{ \sum_{k=1}^n (k-1)/m \cdot s_k(G) \cdot k}{ \sum_{k=1}^n (k-1)/m \cdot s_k(G)} = \frac{ \sum_{k=1}^n k(k-1)\cdot s_k(G)}{ \sum_{k=1}^n (k-1) \cdot s_k(G)}.
    \end{align*}
    Noticing that $$\left( \sum_{k=1}^n k\cdot s_k(G) \right)  \mu(G,v) =\left(\sum_{k=1}^n s_k(G) \right) \mu(G) + \left( \sum_{k=1}^n (k-1) \cdot s_k(G) \right)\mu(G,e),$$
    i.e., that $\mu(G,v)$ is a convex combination (with positive coefficients) of $\mu(G) $ and $\mu(G,e)$, the inequality $\mu(G,v)>\mu(G)$ (\cref{cor:VT_muv>mu}) also implies $\mu(G,e)> \mu(G,v).$
\end{proof}

%\section{Conclusion}\label{sec:conc}

\subsection*{AI declaration}

ChatGPT 5.6 Sol has been used during a revision of this paper, by providing an initial proof of~\cref{conj1}.

% \appendix

% \section{Appendix}

% \begin{lem}\label{lem:VT_muv>mu}
%     In a vertex-transitive graph $G,$ $\mu(G,v)> \mu(G)$ for every $v \in V(G).$
% \end{lem}

% \begin{proof}
%     Every subtree $T_i$ of order $i$ is taken into account in the local subtree order of its $i$ vertices. 
%     Since $v$ is vertex-transitive, $\mu(G,v)$ is constant over all vertices.
%     Hence 
%     $$ \mu(G,v) = \frac{\sum_{i=1}^n i^2 s_i(G)}{\sum_{i=1}^n i s_i(G)} = \frac{ \sum_{k=1}^n \sum_{i=k}^n i s_i(G)}{\sum_{k=1}^n \sum_{i=k}^n s_i(G)} \ge \frac{\sum_{i=1}^n i s_i(G)}{\sum_{i=1}^n s_i(G)} = \mu(G).$$
%     The latter inequality holds since $\frac{ \sum_{i=k}^n i s_i(G)}{\sum_{i=k}^n s_i(G)} $ is an increasing function in $k$ and so the minimum is attained when $i=1$ (which is then smaller than a weighted average with positive coefficients).
% \end{proof}

\end{document}